\documentclass[12pt]{amsart}

\usepackage{amssymb, amsthm, amsmath}

\usepackage{gksymb}

\usepackage[single]{accents}

\newtheorem{theorem}{Theorem}[section]
\newtheorem{lemma}[theorem]{Lemma}

\newtheorem{corollary}[theorem]{Corollary}
\newtheorem{propose}[theorem]{Proposition}

\theoremstyle{definition}
\newtheorem{definition}[theorem]{Definition}

\theoremstyle{remark}
\newtheorem{remark}[theorem]{Remark}

\numberwithin{equation}{section}

\newcommand{\kvec}[2]{\tbinom{#1}{#2}}
\newcommand{\kmat}[4]{\left(\kvec{#1}{#2},\kvec{#3}{#4}\right)}
\newcommand{\refl}[1]{\accentset{\leftarrow}{#1}}
\newcommand{\drefl}[1]{\overleftrightarrow{#1}}
\newcommand{\ip}[2]{\langle #1, #2 \rangle}
\newcommand{\kveczw}{\kvec{z}{w}}
\newcommand{\kmatzw}{\kmat{z}{w}{\zeta}{\omega}}
\newcommand{\gkmeas}{\frac{dz}{z}\frac{dw}{w}}
\newcommand{\Kr}{K_\rho}
\renewcommand{\zeta}{Z}
\renewcommand{\omega}{W}

\title{Bernstein-Szeg\H{o} measures on the two dimensional torus}
\author{Greg Knese} 
\address{Department of Mathematics, Washington
University in St. Louis, St. Louis, MO, 63130}
\date{May 10, 2007} 
\email{geknese@wustl.edu}
\urladdr{http://www.math.wustl.edu/\textasciitilde geknese}
\thanks{Research supported by the Washington University Dissertation
  Fellowship and the Judith Ross Arts \& Sciences Scholarship}

\begin{document}

\begin{abstract}
We present a new viewpoint (namely, reproducing kernels) and new
proofs for several recent results of J. Geronimo and H. Woerdeman on
orthogonal polynomials on the two dimensional torus (and related
subjects).  In addition, we show how their results give a new proof of
And\^{o}'s inequality via an equivalent version proven by Cole and
Wermer.  A major theme is the use of so-called Bernstein-Szeg\H{o}
measures.  A simple necessary and sufficient condition for two
variable polynomial stability is also given.
\end{abstract}

\maketitle

\begin{section}{Introduction}

\begin{subsection}{Prelude}

In several recent papers, J. Geronimo and H. Woerdeman have presented
a number of important generalizations of classical theorems revolving
around such topics as orthogonal polynomials on the unit circle,
Fej\'er-Riesz factorization, and autoregressive filter design to the
context of two variables; namely, generalizations to the two
dimensional torus or just the \emph{2-torus} $\mathbb{T}^2 = (\partial
\mathbb{D})^2 \subset \mathbb{C}^2$ (also known as the distinguished
boundary of the bidisk $\mathbb{D}^2$).  See \cite{Ger-Woe},
\cite{Ger-Woe2}, \cite{Ger-Woe3}, \cite{Ger-Woe4}.  The approach of
Geronimo and Woerdeman may best be described as matrix polynomial
theory.  In this paper we present a new viewpoint---namely that of
reproducing kernels---so that a larger may appreciate what they have
done.

A major theme in this paper is the use of Bernstein-Szeg\H{o}
measures, which are nothing more than absolutely continuous measures
on the 2-torus (or even the $n$-torus for that matter) with weight
given by the reciprocal of the squared modulus of a stable polynomial.
Studying such measures allows us to give a new proof of And\^{o}'s
inequality from operator theory via an equivalent version proven by
B. Cole and J. Wermer.  The proof here is somewhat more elementary
than other proofs (although not necessarily shorter).  As is
well-known, And\^{o}'s inequality is intimately related to finite
interpolation problems for bounded analytic functions on the bidisk
(Pick interpolation), and therefore this work provides a new
connection between ideas in orthogonal polynomials and interpolation.

All of this may sound quite specialized; however, the results in this
paper have consequences that anyone can appreciate.  Let us whet the
reader's appetite with the following interesting by-product of our
work (which will also serve to introduce the terms \emph{degree}
$(n,m)$, \emph{reflection}, and \emph{stable}).

\begin{theorem}\label{stablethm}
Let $q$ be a polynomial in two complex variables of \emph{degree}
$(n,m)$, i.e.
\[
q\kveczw = \sum_{j=0}^{n} \sum_{k=0}^{m} a_{j,k} z^j w^k
\]
and define the \emph{reflection} $\refl{q}$ of $q$ to be 
\[
\refl{q}\kveczw = z^n w^m \overline{q\kvec{1/\bar{z}}{1/\bar{w}}} =
\sum_{j=0}^{n} \sum_{k=0}^{m} \bar{a}_{(n-j),(m-k)} z^j w^k
\]

The polynomial $q$ has no zeros in the closed bidisk
$\overline{\mathbb{D}}^2$ (i.e. $q$ is \emph{stable}) if and only if
there exists a $c>0$ such that
\begin{equation}\label{zerocond}
|q\kveczw|^2 - |\refl{q}\kveczw|^2 \geq c(1-|z|^2)(1-|w|^2)
\end{equation}
for all $z,w \in \overline{\mathbb{D}}$, in which case we can take
\[
\frac{1}{c} = \frac{1}{(2\pi i)^2}\int_{\mathbb{T}^2}
\frac{1}{|q\kveczw|^2}\gkmeas.
\]
\end{theorem}

See Section \ref{sec:Con} for a proof.  Notice that \eqref{zerocond}
certainly implies $q$ has no zeros on the \emph{open} bidisk.  The
point here is that the condition \eqref{zerocond} is exactly what is
needed to imply $q$ has no zeros on the \emph{boundary} of the bidisk.
In some sense, $(1-|z|^2)(1-|w|^2)$ describes the highest possible
order of vanishing of $|q|^2-|\refl{q}|^2$ when $q$ is stable.  At the
same time this says something very specific about two variables:
obvious analogues of this theorem in three variables are probably not
true, but, while the analogue of this theorem in one variable
\emph{is} true, the condition \eqref{zerocond} can be weakened.
Namely, either of the following two weaker conditions are necessary
and sufficient for a one variable polynomial $q$ to have no zeros on
the closed disk:
\[
|q(z)|^2 \geq c(1-|z|^2) \qquad \text{ or } \qquad |q(z)|^2 -
 |\refl{q}(z)|^2 \geq c (1-|z|^2)^2
\]
For the first condition see \cite{Simon}, page 103, and for the second
condition Lemma \ref{lemma:zeros} can be adapted to one variable.

Interestingly enough, the one variable version of this theorem can be
used to prove the fundamental theorem of algebra.  This two variable
version can be used to prove the following: any two variable
polynomial with no zeros on $\overline{\mathbb{D}}^2 \cup
(\mathbb{C}\setminus \mathbb{D})^2$ must be constant.

\end{subsection}

\begin{subsection}{A key to notations and conventions}
All of the following notations will be introduced throughout the
paper, but we collect them here for the benefit of the reader.

Notations/Definitions:
\[
\begin{aligned}
\mathbb{C} &= \text{ complex plane}\\
\mathbb{D} &= \text{ open unit disk in } \mathbb{C}\\
\overline{\mathbb{D}} &= \text{closed unit disk}\\
 \mathbb{T} &=
\partial \mathbb{D} = \text{boundary of the disk in } \mathbb{C}\\
\mathbb{D}^2 &= \mathbb{D} \times \mathbb{D} = \text{ the bidisk}\\
\overline{\mathbb{D}}^2 &= \text{ the closed bidisk}\\ 
\mathbb{T}^2 &=
\mathbb{T} \times \mathbb{T} = \text{ the 2-torus }\\ 
\refl{q}\kveczw &= z^n w^m \overline{q\kvec{1/\bar{z}}{1/\bar{w}}} =
\text{ reflection of } q\\
V \vee W &= \text{ the join or span of two vector spaces}\\
(f\cdot V) &= \{fP: P \in V\} \text{ for a vector space } V\\
q \text{ is stable}& = q \text{ has no zeros in }
\overline{\mathbb{D}}^2\\
\gkbox &=
\text{span}\{z^j w^j: 0\leq j \leq n, 0\leq k\leq m\}\\
 &= \text{
polynomials of degree } (n,m)\\ 
\Kr \gkbox &= \text{ reproducing
kernel for } \gkbox \text{ with respect to a measure } \rho
\end{aligned}
\]

These last two somewhat strange looking notations will make more sense
later.  Many other important notations are defined in Section
\ref{sec:note} but we do not present them here.

Conventions:
\[
\begin{aligned}
& n,m \text{ are positive integers, fixed throughout the paper}\\
& p \text{ is a holomorphic polynomial in two variables of particular importance}\\
& q, P, Q \text{ are generic holomorphic polynomials of two variables}\\
& f \text{ is a typical element of } H^2 \text{ or } L^2\\
& z,w \text{ are the typical holomorphic variables}\\
& \zeta, \omega \text{ are the typical anti-holomorphic variables}\\
& \rho \text{ is a probability measure on } \mathbb{T}^2\\
& \mu \text{ is a Bernstein-Szeg\H{o} measure on } \mathbb{T}^2 \text{
    (see \eqref{bsmeas})}\\
& \kveczw \text{ is an element of } \mathbb{C}^2 \text{ and \emph{not}
    a binomial coefficient}\\
& \frac{1}{2\pi i}\frac{dz}{z} \text{ is normalized Lebesgue measure on } \mathbb{T}\\
& \ip{}{} \text{ is an inner product }\\
& \ip{}{}_\rho \text{ is an inner product given by } \rho\\
& \perp_{\rho} \text{ is orthogonality with respect to } \rho
\end{aligned}
\]

\end{subsection}

\begin{subsection}{Probability measures on $\mathbb{T}^2$}
The story begins with a probability measure $\rho$ on the 2-torus
$(\partial \mathbb{D})^2=\mathbb{T}^2$.  Let $\ip{f}{g}_\rho$ denote
the standard inner product on $L^2(\rho)$ given by
$\int_{\mathbb{T}^2} f\bar{g} d\rho$ and let $||f||$ as usual denote
$\sqrt{\ip{f}{f}_\rho}$.  Fix positive integers $m,n$ throughout the paper,
and define the complex polynomials of degree $(n,m)$
\[
\gkbox := \text{span}\{z^j w^k: 0\leq j \leq n,
0\leq k \leq m\}.
\]

The following nondegeneracy condition will often be imposed on $\rho$:

\begin{definition}\label{ass} The measure $\rho$ is said to be
  \emph{nondegenerate} (at the $(n,m)$ level) if a polynomial $P \in \gkbox$
  can have norm $||P|| =0$ only if $P \equiv 0$.
\end{definition}

This is just another way of saying that the Toeplitz moment matrix for
$\rho$ corresponding to the product index set $\{0,\dots, n\}\times
\{0,\dots, m\}$ is positive definite.  

Let us now define a particular kind of measure that will be of the
utmost importance in this paper.  Recall that a polynomial in two
variables is \emph{stable} if it has no zeros on the closed bidisk.

\begin{definition}\label{def:BS} Given a stable polynomial $q$ the
  \emph{Bernstein-Szeg\H{o}} measure corresponding to $q$ is the
  following probability measure on $\mathbb{T}^2$:
\begin{equation}\label{eq:approx}
d\mu := \frac{c^2}{(2\pi i)^2|q\kveczw|^2} \gkmeas
\end{equation}
where $c>0$ is chosen to make $\mu$ a bona fide probability measure.
\end{definition}

To see why such measures are of interest, let us present a classical
one-variable theorem found in \cite{Simon} (page 95) for instance.

\begin{theorem}[Bernstein-Szeg\H{o} Approximation]
Let $d\nu$ be a nontrivial probability measure (meaning, it is not a
  finite number of point masses) on $\partial \mathbb{D}$ and let $p$
  be the polynomial (in one complex variable $z$) of degree $n$
  satisfying
\[
\begin{aligned}
p &\perp_\nu z , z^2, \dots, z^n\\
p(0) &> 0\\
||p||_\nu &= 1.
\end{aligned}
\]
Then, $p$ has no zeros in the closed disk $\overline{\mathbb{D}}$ and
the measure on $\mathbb{T}$
\[
d\nu_n := \frac{1}{2\pi i |p(z)|^2} \frac{dz}{z}
\]
defines the same inner product on polynomials of degree less than or
equal to $n$ as $d\nu$.

Moreover, $d\nu_n \to d\nu$ in the weak* topology.
\end{theorem}

Obvious analogues of the above theorem do not hold in two variables.
The following theorem of Geronimo and Woerdeman says exactly when a
measure on $\mathbb{T}^2$ has a ``Bernstein-Szeg\H{o}'' approximation.
(Actually this is a slight weakening of their Theorem 1.1.2 in
\cite{Ger-Woe}, but will suffice for our purposes. In Section
\ref{sec:equiv}, we state a group of equivalences that help translate
between our paper and \cite{Ger-Woe}.)

\begin{theorem}[Geronimo-Woerdeman]\label{thm:Ger-Woe}
Suppose $\rho$ is nondegenerate (as in Definition \ref{ass}) and
consider the polynomial $p \in \gkbox$ defined by the conditions
\[
\begin{aligned}
p &\perp_\rho z^j w^k \text{ for } 0\leq j \leq n, 0\leq k \leq m, (i,j)\ne
(0,0)\\
p\kvec{0}{0} &> 0\\
||p|| &= 1.
\end{aligned}
\]
Then, $p$ is stable (i.e. has no zeros in the closed bidisk
$\overline{\mathbb{D}}^2$) and the measure on $\mathbb{T}^2$
\[ 
d\mu := \frac{1}{(2\pi i)^2|p\kveczw|^2} \gkmeas
\] 
defines the same inner product on $\gkbox$ as $\rho$ if and only if
$\rho$ satisfies the following ``automatic orthogonality condition'':
Every polynomial $q$ in the span of $\{z^i w^j: 0\leq i \leq n, 0\leq
j \leq m-1\}$ satisfying
\[
q \perp_\rho z^jw^k \text{ for } 0\leq j\leq n-1, 0\leq k \leq m-1 
\]
is automatically orthogonal to more monomials:
\[
q \perp_\rho z^jw^m \text{ for } 0\leq j\leq n-1. 
\]
\end{theorem}

In Section \ref{sec:note} we shall introduce a new notation and in
Sections \ref{sec:kernel} and \ref{sec:bs} we present several results
which we believe provide insight into this somewhat technical
sounding, yet important, theorem.  Before that, we turn to an
interesting connection between this work, interpolation for bounded
analytic functions on the bidisk, and two-variable operator theory.

\end{subsection}

\begin{subsection}{And\^{o}, Agler, and Christoffel-Darboux}
An important part of the proof of the Geronimo-Woerdeman Theorem
\ref{thm:Ger-Woe} is a ``Christoffel-Darboux like formula'' for two
variables.  To set the stage and for later use, let us present the
Christoffel-Darboux formula for orthogonal polynomials on the unit
circle (see \cite{Simon} Theorem 2.2.7 page 124).

\begin{theorem}[Christoffel-Darboux formula] 
Let $\nu$ be a nontrivial probability measure on $\mathbb{T}$, and let
$p$ be the unit norm polynomial (of one complex variable) of degree
$n$ orthogonal to $z,z^2, \dots, z^n$ with $p(0)>0$.  Also, let
$q_0,q_1, \dots, q_{n-1}$ be any orthonormal basis (with
respect to the inner product given by $\nu$) for the subspace of
polynomials of degree at most $n-1$.  Then,
\[
p(z)\overline{p(\zeta)} - \refl{p}(z)\overline{\refl{p}(\zeta)} =
(1-z\bar{\zeta}) \sum_{j=0}^{n-1} q_j(z) \overline{q_j(\zeta)}
\]
where $\refl{p}(z) := z^n \overline{p(1/\bar{z})}$.
\end{theorem}

The two variable formula proven by Geronimo-Woerdeman has the
following flavor.  For $p$ and $\rho$ as in the conclusion of Theorem
\ref{thm:Ger-Woe}, there exist polynomials $P_j \in \gkbox$,
$j=0,1,\dots n-1$, $Q_k \in \gkbox$, $k=0,\dots, m-1$ such that
\[
\begin{aligned}
p\kvec{z}{w} \overline{p\kvec{\zeta}{\omega}} - \refl{p}\kvec{z}{w}
\overline{\refl{p}\kvec{\zeta}{\omega}} =& (1-z\bar{\zeta})\sum_{j=0}^{n-1}
P_j\kveczw \overline{P_j\kvec{\zeta}{\omega}}\\
&+ (1-w\bar{\omega})
\sum_{k=0}^{m-1} Q_k\kveczw \overline{Q_k\kvec{\zeta}{\omega}}
\end{aligned}
\]
where $\refl{p}\kveczw = z^n w^m \overline{p\kvec{1/\bar{z}}{1/\bar{w}}}$.

For the moment, we do not need to say what $P_j$ and $Q_k$ are
(although we can and will later).  It turns out this
Christoffel-Darboux formula is closely related to And\^{o}'s theorem
in operator theory and Agler's theorem on finite interpolation for
bounded analytic functions on the bidisk.

And\^{o}'s theorem \cite{Ando} says given a polynomial $\sum_{j,k}
a_{jk} z^j w^k$ in two complex variables and two commuting operators
$S,T$ on a Hilbert space satisfying $||S||,||T||\leq 1$, the following
holds:
\[
||\sum_{j,k} a_{jk} S^j T^k|| \leq \sup_{z,w \in \mathbb{T}}
|\sum_{j,k} a_{jk} z^j w^k|.
\]
The Pick interpolation theorem on the bidisk, namely Agler's theorem
(see \cite{Pick} page 180), says given $N$ distinct points
$\lambda_1,\lambda_2,\dots, \lambda_N \in \mathbb{D}^2$ (where we
write $\lambda_i = (z_i,w_i)$) and $N$ points $c_1, \dots, c_n \in
\mathbb{C}$, there exists a holomorphic function $f$ on $\mathbb{D}^2$
with $\sup_{\mathbb{D}^2} |f| \leq 1$ which interpolates
$f(\lambda_i)=c_i$ if and only if there exist positive semidefinite
$N\times N$ matrices $(A_{jk})$ and $(B_{jk})$ so that
\[
1-c_j \bar{c}_k = (1-z_j \bar{z}_k)A_{jk} + (1-w_j \bar{w}_k) B_{jk}
\]
for $1\leq j,k \leq N$.

It turns out that both And\^{o}'s theorem and Agler's theorem are
equivalent to a result about polynomials that looks quite similar to
the Geronimo-Woerdeman formula.

\begin{theorem}[Cole-Wermer \cite{C-W}]\label{CWthm}
Let $P$ and $Q$ be polynomials in two complex variables satisfying
\begin{equation}\label{PgeqQ}
|P\kveczw| \geq |Q\kveczw| \text{ for all } \kveczw \in \mathbb{D}^2
\end{equation}
and
\begin{equation}\label{PeqQ}
|P\kveczw| = |Q\kveczw| \text{ for all } \kveczw \in \mathbb{T}^2.
\end{equation}
Then, there exist polynomials $A_j, B_j$ in two variables, $j=1,\dots,
N$, such that we have
\begin{equation}\label{decomp}
|P\kveczw|^2 - |Q\kveczw|^2 = (1-|z|^2)\sum_{j=1}^{N} |A_j\kveczw|^2 +
 (1-|w|^2)\sum_{j=1}^{N} |B_j\kveczw|^2
\end{equation}
for all $\kveczw \in \mathbb{C}^2$.
\end{theorem}

The conditions \eqref{PgeqQ} and \eqref{PeqQ} are just another way of
saying that $Q/P$ is a rational inner function on the bidisk.
Rational inner functions on the bidisk have a quite simple form: $Q$
as above has to equal $\refl{P}$ (where we have to perform the
``reflection'' at the right degree).  See \cite{Rudin} Theorem
5.2.5. In Section \ref{sec:Ando} we shall present the simple proof
that the Geronimo-Woerdeman formula actually proves the above theorem.
More importantly, the approach of Geronimo-Woerdeman and of this paper
actually gives more information about the decomposition in
\eqref{decomp} than previous proofs of And\^{o}'s theorem. The
decomposition \eqref{decomp} is usually proven with a finite
dimensional Hahn-Banach theorem and as such is not explicit about what
$A_j$ and $B_j$ are. The approach here provides very specific
information about $A_j$ and $B_j$ in \eqref{decomp} (at least in the
case where $P$ is stable).  The key to our proof is the study of
Bernstein-Szeg\H{o} measures.

\end{subsection}
\end{section}

\begin{section}{Notation and Theorem Rephrasing}\label{sec:note}

This section is devoted to a new notation which we believe will make
many of the earlier results more conceptually clear and easier to
digest.  Like earlier, we shall fix a nondegenerate probability
measure $\rho$ on $\mathbb{T}^2$ and corresponding inner product
$\ip{}{}$. Let $V$ be a finite dimensional subspace of two variable
polynomials.  For each $\kvec{\zeta}{\omega} \in \mathbb{C}^2$,
evaluation at $\kvec{\zeta}{\omega}$ is a bounded linear functional on
$V$ (bounded, of course, because $V$ is finite dimensional), and can
therefore be represented as the inner product against the reproducing
kernel for $V$ which we denote as $\Kr V_{\kvec{\zeta}{\omega}}$.
Specifically, if $P\in V$, then $\langle P , \Kr
V_{\binom{\zeta}{\omega}}\rangle = P\kvec{\zeta}{\omega}$.  We define
$\Kr V: \mathbb{C}^2 \times \mathbb{C}^2 \to \mathbb{C}$ without the
subscript $\kvec{\zeta}{\omega}$ to be
\[
\Kr V \kmat{z}{w}{\zeta}{\omega} := \Kr V_{\kvec{\zeta}{\omega}} \kvec{z}{w} =
\langle \Kr V_{\kvec{\zeta}{\omega}} , \Kr V_{\kvec{z}{w}} \rangle,
\]
which happens to be a holomorphic polynomial in $z,w$, an
anti-holo\-morphic polynomial in $\zeta, \omega$, and conjugate
symmetric in the first and second vectors.  Essentially we have
defined a map depending on $\rho$
\[
\begin{aligned}
\Kr:\{\text{finite dimensional subspaces of }\gkbox\} &\to \{\text{
  reproducing kernels}\}\\
 V &\mapsto \text{ reproducing kernel for } V
\end{aligned}
\]
which, as we shall see in Section \ref{sec:kernel}, turns orthogonal
direct sums into sums, vector space shifts into multiplication shifts,
and ``reflection'' of subspaces into ``double reflection'' of
reproducing kernels.  

The following subspaces will be useful throughout.  Once the
definitions are digested, the use of these strange symbols will become
apparent.  These squares should be viewed as lying in a grid in
$\mathbb{Z}^2$ where the lower left corner is the origin and the upper
right corner is the point $(n,m)$.  Essentially we are defining
subspaces of polynomials based on the support of their Fourier
coefficients (and the symbols lie in the Fourier domain
$\mathbb{Z}^2$).

\begin{itemize}
\item Let $\gkboxur$ denote the span of $\{z^j w^k: 0\leq j \leq n,
  0\leq k \leq m, (j,k) \ne (n,m)\}$

\item Let $\gkboxll$ denote the span of $\{z^j w^k: 0\leq j \leq n,
  0\leq k \leq m, (j,k) \ne (0,0)\}$.

\item Let $\gkboxu$ denote the span of $\{z^j w^k : 0\leq j \leq n, 0
  \leq k \leq m-1\}$.

\item Let $\gkboxr$ denote the span of $\{z^j w^k : 0\leq j \leq n-1,
  0\leq k \leq m\}$.

\item Let $\gkboxsm$ denote the span of $\{z^j w^k: 0\leq j \leq n-1,
  0\leq k \leq m-1\}$.

\end{itemize}

We hope now the definitions of the subspaces $\gkboxl$ and $\gkboxd$
  will be clear.  Next, we define many more subspaces in terms of
  orthogonal complements using the inner product $\langle , \rangle$.

Define
\[
\begin{aligned}
\gkboxperprt &:= \gkbox \ominus \gkboxl\\
\gkboxperplt &:= \gkbox \ominus \gkboxr\\
\gkboxurperplt &:= \gkboxur \ominus \gkboxr\\
\gkboxurperp &:= \gkbox \ominus \gkboxur
\end{aligned}
\]
and again we hope the definitions of many other such symbols will be
intuitively clear (in fact, that is the prime reason for this new
notation).  When there is more than one measure present, we will add
an additional subscript to make it clear which measure we are
referring to; e.g. $\gkboxperprt_{\rho}$ refers to the subspace
defined above using the inner product $\ip{}{}_\rho$ defined by
$\rho$.  However, if we are looking at a reproducing kernel $\Kr
\gkboxperprt_{\rho}$ we will leave off the second $\rho$ and just
write $\Kr \gkboxperprt$.

Since the subspace $\gkboxllperp$ is one dimensional, it
is not difficult to see that the corresponding reproducing kernel can
be identified with
\begin{equation}\label{p}
\Kr \gkboxllperp\kmat{z}{w}{\zeta}{\omega} = p\kvec{z}{w}
\overline{p\kvec{\zeta}{\omega}}
\end{equation}
for some unit norm polynomial $p$ (which is unique up to
multiplication by a unimodular constant).

With all of this notation in place let us restate the
Geronimo-Wo\-erdeman Theorem.

\begin{theorem}[Geronimo-Woerdeman]\label{thm:Ger-Woe2}
Given a nondegenerate probability measure $\rho$ on $\mathbb{T}^2$,
any unit norm polynomial $p$ in $\gkboxllperp$ is stable and the
measure on $\mathbb{T}^2$
\begin{equation*}
d\mu := \frac{1}{(2\pi i)^2|p\kveczw|^2} \gkmeas
\end{equation*}
has reproducing kernel $K_\mu \gkbox$ equal to $\Kr \gkbox$ if and
only if
\[
\gkboxurperplt_{\rho} = \gkboxuperplt_{\rho}.
\]
\end{theorem}

\end{section}

\begin{section}{Reproducing Kernel Calculus}\label{sec:kernel}
There are several advantages to working with various reproducing
kernels when studying polynomials with respect to a measure on
$\mathbb{T}^2$.  One advantage is the interface between subspaces and
algebraic operations on kernels demonstrated in the following three
fundamental propositions.  Another advantage is that Bergman identity
(below) provides an interface between orthogonal polynomials and
reproducing kernels.  In this way, choices of ordering polynomials and
orthogonal bases are avoided and absorbed into reproducing kernels.

Like before we fix a nondegenerate probability measure $\rho$ on
$\mathbb{T}^2$ and corresponding inner product $\ip{}{}$.

\begin{propose}\label{prop:sum}  If $V,W$ are finite dimensional subspaces of 
polynomials, then $V$ and $W$ are orthogonal if and only if the join
  of $V$ and $W$, denoted $(V\vee W)$, satisfies the following formula
  for reproducing kernels
\[
\Kr(V\vee W) = \Kr V+\Kr W
\]
\end{propose} 
\begin{proof} Suppose $V$ and $W$ are orthogonal.  Then, for any $P
  \in V\oplus W = V\vee W$, which we may write as $P=Q+R$ for $Q\in V$
  and $R\in W$, we have
\[
\begin{aligned}
\ip{P}{\Kr V_{\kvec{z}{w}}+\Kr W_{\kvec{z}{w}}} &= \ip{Q}{\Kr V_{\kvec{z}{w}}} +
\ip{R}{\Kr W_{\kvec{z}{w}}}\\
& = Q\kvec{z}{w} + R\kvec{z}{w} = P\kvec{z}{w}
\end{aligned}
\]
and therefore, $\Kr(V\oplus W) = \Kr V+ \Kr W$.

For the other direction, first note that the reproducing kernels $\Kr
V_{\kveczw}$ span $V$, since any polynomial in $V$ orthogonal to all
reproducing kernels is identically zero. So, if $\Kr(V\vee W) = \Kr V+\Kr W$ as reproducing
kernels, then $Q \in V$ implies
\[
\ip{Q}{\Kr W_{\kveczw}} = \ip{Q}{\Kr (V\vee W)_{\kveczw} - \Kr V_{\kveczw}} =
Q\kveczw - Q \kveczw = 0
\]
i.e. $Q$ is orthogonal to the reproducing kernel for $W$ and hence all
of $W$.
\end{proof}

For example, $\Kr \gkbox = \Kr(\gkboxr\oplus_\rho \gkboxperplt) =
\Kr \gkboxr + \Kr \gkboxperplt$.

In the following proposition we caution the reader not to confuse the
function $z$ and the variable $z$.

\begin{propose}\label{prop:shift} 
Let $V$ be a finite dimensional subspace of two variable polynomials
and define $(z\cdot V) := \{zP: P\in V\}$ and likewise for $(w\cdot
V)$. Then,
\[
\Kr (z\cdot V)\kmatzw = z\bar{\zeta} \Kr V\kmatzw
\]
and likewise for $(w\cdot V)$.
\end{propose}

\begin{proof} Observe that for any $Q \in V$
\[
\ip{zQ}{ z\bar{\zeta} \Kr V_{\kvec{\zeta}{\omega}}} = \ip{Q}{\bar{\zeta}
  \Kr V_{\kvec{\zeta}{\omega}}} = \zeta Q\kvec{\zeta}{\omega}
\]
since multiplication by $z$ is a unitary.  So, $z\bar{\zeta} \Kr V$
reproduces a generic element $zQ$ of $(z\cdot V)$. The claim follows
by uniqueness of reproducing kernels.
\end{proof}

For example, $\Kr \gkboxl = \Kr(z\cdot \gkboxr) = z\bar{\zeta} \Kr
\gkboxr$ and $\Kr \gkboxdperplt = \Kr (w\cdot \gkboxuperplt) =
w\bar{\omega} \Kr \gkboxuperplt$.

For the third proposition, we recall one definition and introduce
another related one.

\begin{definition} Given any $P\in \gkbox$, the \emph{reflection} of
  $P$ (at the $(n,m)$ level) is another polynomial $\refl{P} \in
  \gkbox$ which is defined by
\[
\refl{P}\kvec{z}{w} := z^n w^m
\overline{P\kvec{1/\bar{z}}{1/\bar{w}}}
\]
\end{definition}

\begin{remark}\label{rem:au} 

The map $P \mapsto \refl{P}$ is anti-unitary; i.e. 
\[
\ip{\refl{P}}{\refl{Q}} = \ip{Q}{P}
\]
and can be described in concrete terms as the anti-linear map
which sends a monomial $z^j w^k \mapsto z^{n-j} w^{m-k}$. 

\end{remark}

\begin{definition} Given $\Delta:\mathbb{C}^2 \times \mathbb{C}^2 \to
  \mathbb{C}$ a holomorphic polynomial in the first two variables of
  degree $(n,m)$ and an anti-holomorphic polynomial in the second two
  variables of degree $(n,m)$, the \emph{double reflection} of
  $\Delta$ shall be denoted $\drefl{\Delta}$ and defined by
\[
\drefl{\Delta} \kmat{z}{w}{\zeta}{\omega} := (z\bar{\zeta})^n (w
\bar{\omega})^m \Delta
\kmat{1/\bar{\zeta}}{1/\bar{\omega}}{1/\bar{z}}{1/\bar{w}} 
\]
\end{definition}

\begin{remark}\label{rem:cs} If $\Delta$ is conjugate symmetric, then 
\[
\drefl{\Delta}\kmatzw = \bar{\zeta}^n \bar{\omega}^m
\overleftarrow{\Delta_{\kvec{1/\bar{\zeta}}{1/\bar{\omega}}}}
\kveczw
\]
where $\Delta_{\kvec{1/\bar{\zeta}}{1/\bar{\omega}}} \kveczw = \Delta
\kmat{z}{w}{1/\bar{\zeta}}{1/\bar{\omega}}$.
\end{remark}

\begin{propose}\label{refl:prop} Let $V$ be a subspace of $\gkbox$ and define $\refl{V}
 := \{\refl{P} : P\in V\}$.  Then,
\[
\Kr \refl{V} = \drefl{\Kr V} 
\]
\end{propose}

\begin{proof}
Let $P \in V$.  Then, $\refl{P} \in \refl{V}$ and by Remarks
\ref{rem:au} and \ref{rem:cs}
\[
\ip{\refl{P}}{\drefl{\Kr V}_{\kvec{\zeta}{\omega}}} =
\ip{\refl{P}}{\bar{\zeta}^n \bar{\omega}^m
\overleftarrow{\Kr V_{\kvec{1/\bar{\zeta}}{1/\bar{\omega}}
}}} = \zeta^n \omega^m \ip{\Kr
V_{\kvec{1/\bar{\zeta}}{1/\bar{\omega}}}}{P} =
\refl{P}\kvec{\zeta}{\omega}
\]
and the claim follows, again by uniqueness.
\end{proof}

\begin{propose}[Bergman Identity]\label{prop:Berg} Let $V$ be a finite dimensional
  subspace of two variable polynomials and let $q_1, q_2,\dots, q_N$
  be an orthonormal basis for $V$.  Then,
\[
\Kr V\kmatzw = \sum_{j=1}^{N} q_j\kveczw
\overline{q_j\kvec{\zeta}{\omega}}
\]
\end{propose}

\begin{proof} The expression on the right reproduces each orthonormal
  basis element and hence everything in $V$.
\end{proof}

With the fundamental properties of these reproducing kernels out of
the way let us present a formula that holds in complete generality and
sheds some light on the Geronimo-Woerdeman theorem.

Since $\gkboxllperp$ is one dimensional, the Bergman identity tells us
there is a unit norm polynomial $p$, unique up to multiplication by a
unimodular constant, such that
\[
\Kr \gkboxllperp \kmat{z}{w}{\zeta}{\omega} = p\kvec{z}{w}
\overline{p\kvec{\zeta}{\omega}}.
\]
Allow us to abbreviate this as $\Kr \gkboxllperp = p\bar{p}$.  By
Proposition \ref{refl:prop} above $\Kr \gkboxurperp = \refl{p}
\bar{\refl{p}}$.

\begin{theorem}\label{fundthm} Let $\rho$ be a nondegenerate
  probability measure on $\mathbb{T}^2$, and let $p$ be any unit norm
  polynomial in $\gkboxllperp$.  Then, suppressing the argument
  $\kmatzw$ we have
\[
\begin{aligned}
p\bar{p} - \refl{p}\bar{\refl{p}} =& (1-z\bar{\zeta})\Kr \gkboxrperpdn
+ (1-w\bar{\omega}) \Kr \gkboxuperplt +
(1-z\bar{\zeta})(1-w\bar{\omega})\Kr \gkboxsm \\ &+ \left[\left(\Kr
\gkboxurperplt - \Kr \gkboxuperplt\right) - \left(\Kr \gkboxllperprt -
\Kr \gkboxdperprt\right)\right].
\end{aligned}
\]
\end{theorem}

\begin{proof} Throughout the proof we will suppress $\Kr$ and the
  argument $\kmatzw$.  Observe that
\[
\begin{aligned}
\gkboxperprt - \gkboxperplt &= 
\left( \gkbox - \gkboxl\right) - \left(\gkbox - \gkboxr\right)\\
&= \gkboxr - \gkboxl.
\end{aligned}
\]
But, $\gkboxl = z\bar{\zeta} \gkboxr$ by Proposition \ref{prop:shift}
and therefore,
\begin{equation}\label{fundthm:eq0}
\gkboxperprt - \gkboxperplt = (1-z\bar{\zeta}) \gkboxr
\end{equation}
and similarly
\begin{equation}\label{fundthm:eq1}
\gkboxuperprt - \gkboxuperplt = (1-z\bar{\zeta})\gkboxsm.
\end{equation}
Subtracting the right hand side of equation \eqref{fundthm:eq1} from
the right hand side of equation \eqref{fundthm:eq0} yields
\[
(1-z\bar{\zeta}) (\gkboxr - \gkboxsm) = (1-z\bar{\zeta}) \gkboxrperpdn
\]
by Proposition \ref{prop:sum}.

Subtracting the left hand side of equation \eqref{fundthm:eq1} from
the left hand side of equation \eqref{fundthm:eq0} yields
\[
\begin{aligned}
\left(\gkboxperprt -
\gkboxperplt\right) - \left(\gkboxuperprt - \gkboxuperplt\right) &=
\left(\gkboxperprt - \gkboxdperprt\right) +
\left(\gkboxdperprt - \gkboxuperprt\right)
- \left(\gkboxperplt - \gkboxuperplt\right)\\ 
&=\left(\gkboxperprt - \gkboxdperprt\right) 
- (1-w\bar{\omega})\gkboxuperprt -
\left(\gkboxperplt - \gkboxuperplt\right)
\end{aligned}
\]
where the last equality comes from $w\bar{\omega} \gkboxuperprt =
\gkboxdperprt$ (i.e. Proposition \ref{prop:shift}).

So the difference of equations \eqref{fundthm:eq0} and
\eqref{fundthm:eq1} is 
\[
(1-z\bar{\zeta}) \gkboxrperpdn = \left(\gkboxperprt - \gkboxdperprt\right) 
- (1-w\bar{\omega})\gkboxuperprt -
\left(\gkboxperplt - \gkboxuperplt\right).
\]

Therefore,
\begin{align}
\label{fundthm:eq2}
\left(\gkboxperprt - \gkboxdperprt\right)-\left(\gkboxperplt -
\gkboxuperplt\right) &= (1-z\bar{\zeta})\gkboxrperpdn &+&
(1-w\bar{\omega})\gkboxuperprt\\ 
\nonumber &=
(1-z\bar{\zeta})\gkboxrperpdn 
&+& (1-w\bar{\omega})\gkboxuperplt\\
\nonumber & &+& (1-z\bar{\zeta})(1-w\bar{\omega})\gkboxsm
\end{align}

where the last equality follows from \eqref{fundthm:eq1}.
Finally, if we note that
\[
\gkboxperprt = \gkboxllperp + \gkboxllperprt = p\bar{p} +
\gkboxllperprt
\]
and
\[
\gkboxperplt = \refl{p}\bar{\refl{p}} + \gkboxurperplt
\]
the theorem follows from \eqref{fundthm:eq2} after a little
rearrangement.  Namely, write
\[
\left(\gkboxperprt - \gkboxdperprt\right)-\left(\gkboxperplt -
\gkboxuperplt\right) = p\bar{p} - \refl{p} \bar{\refl{p}} -\left[
\left(\gkboxurperplt - \gkboxuperplt \right) -\left(\gkboxllperprt -
\gkboxdperprt\right) \right]
\]
and use \eqref{fundthm:eq2} to obtain the theorem.
\end{proof}

This easily yields the following corollary.

\begin{corollary}\label{cor:ineq}
With the same setup as Theorem \ref{fundthm}, if $\gkboxurperplt =
\gkboxuperplt$, then
\[
|p\kveczw|^2 - |\refl{p}\kveczw|^2 \geq (1-|z|^2)(1-|w|^2)
\]
\end{corollary}
\begin{proof} When $\gkboxurperplt = \gkboxuperplt$, it is also true that  
$\gkboxllperprt = \gkboxdperprt$ by reflection.  Therefore, on the
diagonal $z=\zeta$, $w=\omega$
\[
|p\kveczw|^2 - |\refl{p} \kveczw|^2 \geq (1-|z|^2)(1-|w|^2)\Kr
 \gkboxsm. 
\]
If we use $\mathbb{C}$ to denote the subspace of constant polynomials,
then the kernel $\Kr \gkboxsm$ can be split into $\Kr \mathbb{C} + \Kr
(\gkboxsm \ominus \mathbb{C})$.  Since $\rho$ is a probability
measure, $\Kr \mathbb{C} \equiv 1$, and therefore, on the diagonal
$z=\zeta$, $w=\omega$, $\Kr\gkboxsm \geq 1$.  Hence,
\[
|p\kveczw|^2 - |\refl{p} \kveczw|^2 \geq (1-|z|^2)(1-|w|^2).
\]
\end{proof}

This corollary combines with the following lemma about polynomials in
two variables to yield one part of the Geronimo-Woerdeman theorem.
This lemma is one direction of the theorem from the introduction.

\begin{lemma}\label{lemma:zeros}
Let $q \in \gkbox$ and suppose there is a $c>0$ such that
\[
|q\kveczw|^2 - |\refl{q}\kveczw|^2 \geq c(1-|z|^2)(1-|w|^2)
\]
for all $z,w \in \mathbb{D}$. Then, $q$ is stable.
\end{lemma}


\begin{proof}
First, suppose $q$ has a zero, say $\kvec{z_0}{w_0}$, on $\mathbb{T}
\times \mathbb{D}$. Then, $|q\kvec{rz_0}{w_0}|^2 = O(1-r)^2$ but
\[
|q\kvec{rz_0}{w_0}|^2 \geq c(1-|w_0|^2)(1-r^2)
\]
which fails as $r\nearrow 1$.  Similarly, $q$ has no zeros on
$\mathbb{D}\times \mathbb{T}$.

Second, suppose $q$ has a zero, $\kvec{z_0}{w_0}$, on $\mathbb{T}^2$.
Writing $q\kvec{rz_0}{rw_0} = a(1-r) + O(1-r)^2$ we have
\[
\refl{q}\kvec{rz_0}{rw_0} = r^{n+m} z_{0}^n w_{0}^m
\overline{q\kvec{z_0/r}{w_0/r}} = 
r^{n+m} z_{0}^n
w_{0}^m\left(\bar{a}(1-\frac{1}{r}) + O(1-r)^2\right).
\]

Therefore,
\[
\begin{aligned}
|q\kvec{rz_0}{rw_0}|^2 - |\refl{q}\kvec{rz_0}{rw_0}|^2 =&
 |a(1-r)+O(1-r)^2|^2\\
 &- r^{2(n+m-1)}|a(1-r)+O(1-r)^2|^2\\
=& |a|^2(1-r)^2(1-r^{2(n+m-1)})+O(1-r)^3\\
=& O(1-r)^3 \geq  c(1-r^2)^2
\end{aligned}
\]
which also fails as $r\nearrow 1$.
\end{proof}

\begin{corollary}
With the setup of Theorem \ref{fundthm}, if $\gkboxurperplt =
\gkboxuperplt$, then $p$ is stable.
\end{corollary}

This proves one part of the reverse implication of Theorem
\ref{thm:Ger-Woe2}.  For the remainder of the reverse implication, we
refer the reader to the original paper \cite{Ger-Woe}.

\end{section}

\begin{section}{Bernstein-Szeg\H{o} Measures on $\mathbb{T}^2$} \label{sec:bs}
Next, we would like to indicate why the forward implication of the
Geronimo-Woerdeman theorem holds.  This involves studying
Bern\-stein-Szeg\H{o} measures:
\begin{equation}\label{bsmeas}
d\mu := \frac{c^2}{(2\pi i)^2|q\kveczw|^2} \gkmeas
\end{equation}
where $q$ is a polynomial in $\gkbox$ with no zeros in the closed
bidisk and $c>0$ is chosen to make $\mu$ a probability measure.

The first thing to notice about such a measure is that the norm it
provides for $L^2(\mathbb{T}^2, \mu)$ is equivalent to the norm on
$L^2(\mathbb{T}^2)$ with Lebesgue measure (because $q$ is bounded
above and below on $\mathbb{T}^2$).  So, all closures taken with
respect to Lebesgue measure are equal to closures taken with respect
to the norm determined by $\mu$.  For instance,
$H^2(\mathbb{T}^2,\mu)$, the closure of the polynomials in
$L^2(\mathbb{T}^2,\mu)$, is equal to the usual Hardy space
$H^2(\mathbb{T}^2)$.

Another important fact about Bernstein-Szeg\H{o} measures is given in
the following proposition.

\begin{propose}\label{bsqprop} For $\mu$ and $q$ as in \eqref{bsmeas},
\[
q \in \gkboxllperp \text{ and } K_\mu \gkboxllperp \kmatzw =
\frac{1}{c^2} q\kveczw \overline{q\kvec{\zeta}{\omega}}
\]
\end{propose}
\begin{proof} Observe that 
\[
\begin{aligned}
\ip{z^jw^k}{q} &= \frac{c^2}{(2\pi i)^2} \int_{\mathbb{T}^2} \frac{z^j
  w^k \overline{q\kveczw}}{|q\kveczw|^2} \gkmeas\\
&=  \frac{c^2}{(2\pi i)^2} \int_{\mathbb{T}^2} \frac{z^j
  w^k }{q\kveczw} \gkmeas
\end{aligned}
\]
and this equals zero when $j>0$ or $k>0$.  So, $q \in \gkboxllperp$
and since $||q|| =c$, it follows that $K_\mu \gkboxllperp \kmatzw =
\frac{1}{c^2} q\kveczw \overline{q\kvec{\zeta}{\omega}}$ by the
Bergman identity.
\end{proof}

Let us define several important closed subspaces of
$H^2(\mathbb{T}^2)$. 

\begin{definition} All closed spans below are taken with respect to
  $H^2(\mathbb{T}^2)$.
\begin{enumerate} 
\item The \emph{border subspace} $\mathbf{B}$ is defined to be
\[
\mathbf{B} := \overline{\text{span}}\{z^j w^k: j, k\geq 0 \text{ and
  either } j\leq n-1 \text{ or } k \leq m-1\}.
\]

\item The \emph{bottom border subspace} $\mathbf{BB}$ is defined to be
\[
\mathbf{BB} := \overline{\text{span}}\{z^j w^k: j\geq 0 \text{ and }
0\leq k \leq m-1\}.
\]

\item The \emph{left border subspace} $\mathbf{LB}$ is defined to be
\[
\mathbf{LB} := \overline{\text{span}}\{z^j w^k: k\geq 0 \text{ and }
0\leq j \leq n-1\}.
\]
\end{enumerate}
\end{definition}

We highly recommend drawing pictures of the above subspaces analogous
to the pictures from earlier.

\begin{lemma}\label{bs:lemma} 
Define 
\[
L_{\kvec{\zeta}{\omega}} \kveczw := \frac{(z\bar{\zeta})^n
  [q\kveczw\overline{q\kvec{1/\bar{z}}{\omega}} - \refl{q}\kveczw
  \overline{\refl{q}\kvec{1/\bar{z}}{\omega}}]}{(1-z\bar{\zeta})(1-w\bar{\omega})}
\]
The functions $L_{\kvec{\zeta}{\omega}}$ for $\zeta,\omega \in
  \mathbb{D}$ have the following properties.

\begin{enumerate}
\item $L_{\kvec{\zeta}{\omega}}$ is orthogonal (with respect to $\mu$)
  to $\mathbf{LB}$.

\item $L_{\kvec{\zeta}{\omega}}$ is an element of $\mathbf{BB}$.

\item $L_{\kvec{\zeta}{\omega}}$ reproduces $(z^n \cdot
  \mathbf{BB})$; i.e. $f\kveczw = \ip{f}{L_{\kveczw}}_\mu$ for $f$ in the
  closed span of $\{z^j w^k: j \geq n, 0\leq k \leq m-1\}$.

\end{enumerate}
\end{lemma}

\begin{remark}\label{rem:ess} The functions $L_{\kveczw}$ are
  purposefully asymmetric: the variable $z$ appears six times in the
  definition and $\zeta$ appears only twice.  In essence, the
  proposition says that $L_{\kveczw}$ reproduces the
  $H^2(\mathbb{T}^2)$ projection of $f \in \mathbf{B}$ to $(z^n \cdot
  \mathbf{BB})$.
\end{remark}

Before proving this let us show how it proves the following important
fact for Bernstein-Szeg\H{o} measures.

\begin{theorem}\label{bs:thm} Let $q$ be a polynomial with no zeros on
  $\overline{\mathbb{D}}^2$ and define its Bernstein-Szeg\H{o} measure
  $\mu$ as in \eqref{bsmeas}.  Then, 
\[
\begin{aligned}
\mathbf{B} \ominus_\mu \mathbf{BB} &= \mathbf{LB} \ominus_\mu
\gkboxsm ,\\
\mathbf{B} \ominus_\mu \mathbf{LB} &= \mathbf{BB} \ominus_\mu
\gkboxsm
\end{aligned}
\]
and
\[
\gkboxurperpdn_\mu = \gkboxrperpdn_\mu \text{ and } \gkboxurperplt_\mu
= \gkboxuperplt_\mu
\]

\end{theorem}

\begin{proof}
Let $f \in \mathbf{B} \ominus_\mu \mathbf{BB}$.  Now, since $f \in
\mathbf{B}$, we can decompose $f = g + h$ where $g$ is the orthogonal
projection of $f$ with respect to \emph{Lebesgue} measure on
$\mathbb{T}^2$ to the subspace $(z^n \cdot \mathbf{BB})$ and $h$ is
the orthogonal projection of $f$ to the subspace $\mathbf{LB}$, again
with \emph{Lebesgue} measure.  So, on the one hand, $f \perp_\mu
\mathbf{BB}$ and $L_{\kveczw} \in \mathbf{BB}$ imply
$\ip{f}{L_{\kveczw}}_\mu = 0$ for all $\kveczw \in \mathbb{D}^2$.  On
the other hand,
\[
\ip{g+h}{L_{\kveczw}}_\mu = g\kveczw
\]
by (1) and (3) of Lemma \ref{bs:lemma}.  Hence, $g \equiv 0$ and
therefore $f = h \in \mathbf{LB}$. So, $\mathbf{B} \ominus_\mu
\mathbf{BB} \subseteq \mathbf{LB} \ominus_\mu \gkboxsm$.  

To prove equality suppose that $f\in \mathbf{LB} \ominus_\mu \gkboxsm$
and $f \perp_\mu \mathbf{B} \ominus_\mu \mathbf{BB}$.  Then, $f \in
\mathbf{BB} \cap \mathbf{LB} = \gkboxsm$, yet at the same time $f
\perp_\mu \gkboxsm$; i.e. $f\equiv 0$.  Hence, $\mathbf{B} \ominus_\mu
\mathbf{BB} = \mathbf{LB} \ominus_\mu \gkboxsm$.  It follows easily
from this that $\gkboxurperpdn_\mu = \gkboxrperpdn_\mu$. The rest of
the theorem follows by symmetry.
\end{proof}

Now, we prove Lemma \ref{bs:lemma} via a series of propositions, which
are interesting in their own right.

\begin{propose}\label{projprop} Define $F_{\kvec{\zeta}{\omega}} \in H^2(\mathbb{T}^2)$ by
\[
F_{\kvec{\zeta}{\omega}} \kveczw = \frac{(z\bar{\zeta})^n
q\kveczw\overline{q\kvec{1/\bar{z}}{\omega}}}{c^2(1-z\bar{\zeta})(1-w\bar{\omega})}
\]
Then, for any $f \in H^2(\mathbb{T}^2)$, the function
\[
g\kveczw := \ip{f}{F_{\kveczw}}_\mu
\]
is the $H^2(\mathbb{T}^2)$-orthogonal projection of $f$ to the
subspace $(z^n \cdot H^2(\mathbb{T}^2))$. (So, a projection in $H^2$
with Lebesgue measure is achieved via the $\mu$-inner product.) 
\end{propose}
\begin{proof} Observe that
\[
\begin{aligned}
g\kvec{\zeta}{\omega} &= \frac{1}{(2\pi i)^2} \int_{\mathbb{T}^2}
  f\kveczw \frac{(\bar{z}\zeta)^n \overline{q\kveczw}q\kvec{z}{
  \omega}}{(1-\bar{z}\zeta)(1-\bar{w}\omega)|q\kveczw|^2}
  \frac{dw}{w}\frac{dz}{z}\\ 
&= \frac{1}{2\pi i}\int_{\mathbb{T}} f\kvec{z}{\omega}
  \frac{(\bar{z}\zeta)^n}{1-\bar{z}\zeta} \frac{dz}{z}\\
\end{aligned}
\]
and the last line does indeed equal the projection of $f$ to the
subspace $(z^n \cdot H^2(\mathbb{T}^2))$.
\end{proof}

\begin{propose}\label{orthprop} The closed subspace $(\refl{q} \cdot
  H^2(\mathbb{T}^2))$ equals $H^2(\mathbb{T}^2) \ominus_\mu
  \mathbf{B}$.
\end{propose}
\begin{proof}  If $f \in H^2(\mathbb{T}^2)$ then
\[
\begin{aligned}
\ip{f \refl{q}}{z^j w^k}_\mu &= \ip{f z^{n-j} w^{m-j}}{q}_\mu \\ 
&= \frac{c}{(2\pi i)^2}\int_{\mathbb{T}^2} f\kveczw
\frac{z^{n-j}w^{m-k}}{q\kveczw} \gkmeas \\
\end{aligned}
\]
equals zero if $j < n$ or $k < m$.  So, $(\refl{q} \cdot
H^2(\mathbb{T}^2)) \subseteq H^2(\mathbb{T}^2) \ominus_\mu
\mathbf{B}$.  On the other hand, if $f \in H^2(\mathbb{T}^2)
\ominus_\mu \mathbf{B}$ and $f \perp_\mu (\refl{q} \cdot
H^2(\mathbb{T}^2))$, then since $\refl{q} \in \gkbox$ and the
coefficient of $z^n w^m$ in $\refl{q}$ is necessarily nonzero, it
follows that $f \perp_\mu z^n w^m$.  But, since $f \perp_\mu
w\refl{q}$, this implies $f \perp_\mu z^n w^{m+1}$.  Continuing like
this we see that $f \perp_\mu z^j w^k$ for all $j,k \geq 0$; i.e. $f
\equiv 0$.  So, $(\refl{q} \cdot H^2(\mathbb{T}^2)) =
H^2(\mathbb{T}^2) \ominus_\mu \mathbf{B}$.
\end{proof}

\begin{propose}\label{polyprop} The function 
\[
G\kmatzw := \frac{q\kveczw \overline{q\kvec{\zeta}{\omega}} -
  (z\bar{\zeta})^n \refl{q}\kvec{1/\bar{\zeta}}{w}
  \overline{\refl{q}\kvec{1/\bar{z}}{\omega}}}{1-w\bar{\omega}}
\]
is a polynomial of degree $m-1$ in $w$ and $\bar{\omega}$.  
\end{propose}
 
\begin{proof} The proposition follows for algebraic reasons, but we
  can actually give meaning to $G$ so we provide the following more
  explicit proof.  

First, notice that if the claim is true when $z=\zeta$ then it is true
  in general.  This is because of the polarization theorem for
  holomorphic functions (i.e. $f(z,\bar{z})\equiv 0$ implies $f \equiv
  0$) applied to coefficients of $w^j\bar{\omega}^k$ in the expansion
  of $G$ above.

To prove the result for $z=\zeta$ we can appeal to the one variable
Christoffel-Darboux formula.  Namely, for each $z \in
\overline{\mathbb{D}}$ let
\[
d\sigma_z(w) := \frac{c^2_z}{2\pi i|q\kveczw|^2} \frac{dw}{w}
\]
where $c_z > 0$ is chosen to make $\sigma_z$ a probability measure on
$\mathbb{T}$.  Then, it is not difficult to check the polynomial in
$w$ given by $q\kveczw/c_z$ satisfies the conditions of the
Christoffel-Darboux theorem (up to multiplication by a unimodular
which does not matter).  Therefore, if we denote the reproducing
kernel for $\sigma_z$ for polynomials in $w$ with degree up to $m-1$
by $K_{m-1}^{z}$, then
\[
\begin{aligned}
c^2_z(1-w\bar{\omega})K_{m-1}^{z} (w,\omega) &=
q\kveczw\overline{q\kvec{z}{\omega}}- (w\bar{\omega})^m
\overline{q\kvec{z}{1/\bar{w}}} q\kvec{z}{1/\bar{\omega}} \\
&= q\kveczw\overline{q\kvec{z}{\omega}} - (z\bar{\zeta})^n \refl{q}\kvec{1/\bar{z}}{
  w} \overline{\refl{q}\kvec{1/\bar{z}}{\omega}}
\end{aligned}
\]
and this proves $G$ is a polynomial in $w$ and $\bar{\omega}$ of
degree $m-1$ for $z=\zeta \in \overline{\mathbb{D}}$.
\end{proof}

These three propositions allow us to prove Lemma \ref{bs:lemma} with
ease.

\begin{proof}[Proof of Lemma \ref{bs:lemma}]
\begin{enumerate}

\item The functions $L_{\kveczw}$ equal a sum of a functions
  orthogonal to $\mathbf{LB}$ by Propositions \ref{projprop} and
  \ref{orthprop}.

\item The function $L_{\kvec{\zeta}{\omega}}$ belongs to $\mathbf{BB}$ by
  Proposition \ref{polyprop} because
\[
L_{\kvec{\zeta}{\omega}} \kveczw = \frac{(z\bar{\zeta})^n}{(1-z\bar{\zeta})}
G\kmat{z}{w}{1/\bar{z}}{\omega}
\]
is holomorphic in $z, \bar{\zeta}, w, \bar{\omega}$ and a polynomial
of degree $m-1$ in $w, \bar{\omega}$.

\item Finally, $L_{\kveczw}$ is the sum of a function that reproduces
  $(z^n \cdot H^2(\mathbb{T}^2))$ and a function orthogonal to $(z^n
  \cdot \mathbf{BB})$.  So, $L_{\kveczw}$ reproduces functions in
  $(z^n \cdot \mathbf{BB})$.
\end{enumerate}

\end{proof}

Finally, Theorem \ref{bs:thm} allows us to prove the forward
implication of Theorem \ref{thm:Ger-Woe2}, which we state as a
corollary.

\begin{corollary} 
Given a nondegenerate probability measure $\rho$ on $\mathbb{T}^2$, if
a unit norm polynomial $p \in \gkboxllperp$ is stable and the measure
on $\mathbb{T}^2$
\begin{equation*}
d\mu := \frac{1}{(2\pi i)^2|p\kveczw|^2} \gkmeas
\end{equation*}
has reproducing kernel $K_\mu \gkbox$ equal to $\Kr \gkbox$, then
\[
\gkboxurperplt_\rho = \gkboxuperplt_\rho.
\]
\end{corollary}
\begin{proof} If $\Kr \gkbox$ equals $K_\mu \gkbox$, then Theorem
  \ref{bs:thm} proves
\[
\gkboxurperplt_{\rho} = \gkboxurperplt_{\mu} = \gkboxuperplt_{\mu} =
\gkboxuperplt_{\rho}.
\]
\end{proof}

\end{section} 

\begin{section}{And\^{o}, Agler, and Christoffel-Darboux Redux}\label{sec:Ando}
Let us present and prove the Ger\-onimo-Woerde\-man
Christ\-offel-Dar\-boux formula (plus another identity) in our own
language.

\begin{theorem} If $\rho$ is a nondegenerate probability measure on
  $\mathbb{T}^2$ and $\gkboxurperplt = \gkboxuperplt$, then writing,
  as usual, $p\bar{p} = \Kr \gkboxllperp$ and suppressing $\kmatzw$
\[
\begin{aligned}
p\bar{p} - \refl{p}\bar{\refl{p}} &= (1-z\bar{\zeta})\Kr \gkboxrperpdn +
(1-w\bar{\omega}) \Kr \gkboxuperplt +
(1-z\bar{\zeta})(1-w\bar{\omega})\Kr \gkboxsm \\
&= (1-z\bar{\zeta})\Kr \gkboxrperpdn +
(1-w\bar{\omega}) \Kr \gkboxuperprt\\
&=(1-z\bar{\zeta})\Kr \gkboxrperpup +
(1-w\bar{\omega}) \Kr \gkboxuperplt.
\end{aligned}
\]
\end{theorem}
\begin{proof} These all follow from Theorem \ref{fundthm} using the
  fact that $\gkboxurperplt = \gkboxuperplt$ implies $\gkboxllperprt =
\gkboxdperprt$ (by reflection) and the
identities 
\[
\begin{aligned}
\Kr \gkboxuperprt - \Kr \gkboxuperplt &= (1-z\bar{\zeta})\Kr \gkboxsm\\
\Kr \gkboxrperpup - \Kr \gkboxrperpdn &= (1-w\bar{\omega})\Kr \gkboxsm.
\end{aligned}
\]
\end{proof}

Let us now prove the Cole-Wermer equivalent of And\^{o}'s theorem.

\begin{proof}[Proof of Theorem \ref{CWthm}] Let us recall the setup.  We have
  two-variable polynomials $P$ and $Q$ satisfying
\[
|P| = |Q| \text{ on } \mathbb{T}^2
\]
and
\[
|P| \geq |Q| \text{ on } \mathbb{D}^2
\]
There is no loss in assuming $P$ and $Q$ are relatively prime (just
multiply the formula we prove by any common factors if necessary).
Then, $Q/P$ is holomorphic on $\mathbb{D}^2$ because it is bounded by
1 wherever it is defined.  In fact, $P$ has no zeros on
$\mathbb{D}^2$, because, first of all, any zeros of $P$ in
$\mathbb{D}^2$ necessarily coincide with zeros of $Q$, since $Q/P$ is
holomorphic.  So, the variety $\{P=0\}\cap\mathbb{D}^2$ (which has
infinitely many points if it is nonempty) is contained in
$\{Q=0\}\cap\mathbb{D}^2$.  But, two relatively prime polynomials in
two variables have only finitely many common zeros (by Bezout's
theorem for instance).  Hence, $\{P=0\}\cap \mathbb{D}^2$ is the empty
set.  This implies $Q/P$ is a rational inner function.  It is a fact
that $Q = \refl{P}$ where we have to perform the ``reflection'' at the
appropriate level, which we assume to be $(n,m)$ (see \cite{Rudin}
Theorem 5.2.5).

If it were the case that $P$ had no zeros on
$\overline{\mathbb{D}}^2$, then the theorem would be proved by
defining a Bernstein-Szeg\H{o} measure using $P$ and applying the
above Christoffel-Darboux formula.  To handle the case where $P$ might
have zeros on $\mathbb{T}^2$ we look at $P_r\kveczw := P\kvec{rz}{rw}$ and let
$r \nearrow 1$.

Define $\mu_r$ to be the Bernstein-Szeg\H{o} measure for $P_r$, i.e.
\[
d\mu_r := \frac{c^2_r}{(2\pi i)^2 |P_r\kveczw)|^2} \gkmeas
\]
where $c_r$ is chosen to make $\mu_r$ a probability measure.  Then, by
Proposition \ref{bsqprop} and the Christoffel-Darboux formula on the
diagonal $z=\zeta$, $w=\omega$
\[
|P_r|^2 - |\refl{P_r}|^2 = c_r^2(1-|z|^2)K_{\mu_r} \gkboxrperpdn +
c_r^2(1-|w|^2) K_{\mu_r} \gkboxuperprt.
\]
Also, by the Bergman identity (evaluating $K_{\mu_r} \gkboxrperpdn$ at
$\kmat{z}{w}{z}{w}$)
\[
\begin{aligned}
\frac{c_r^2}{(2\pi i)^2}\int_{\mathbb{T}^2}K_{\mu_r} \gkboxrperpdn
  \gkmeas 
&\leq 
\left(\sup_{\mathbb{T}^2} |P_r|^2\right)
  \frac{c_r^2}{(2\pi i)^2} \int_{\mathbb{T}^2}
  \frac{K_{\mu_r} \gkboxrperpdn}{|P_r|^2} \gkmeas \\
&\leq
  n\sup_{\mathbb{T}^2} |P|^2
\end{aligned}
\]
since the dimension of $\gkboxrperpdn$ is $n$. Similarly,
\[
\frac{c_r^2}{(2\pi i)^2}\int_{\mathbb{T}^2} K_{\mu_r} \gkboxuperprt
\gkmeas \leq m \sup_{\mathbb{T}^2} |P|^2.
\]

Hence, the polynomials forming $c_r^2 K_{\mu_r} \gkboxrperpdn$ as in the
Bergman identity are bounded in $L^2$ norm and are all of bounded
degree (likewise for $c_r^2 K_{\mu_r} \gkboxuperprt$).  Taking
subsequences if necessary, said polynomials will converge to
polynomials as $r\nearrow 1$.  Also, $\refl{P_r}$ tends to $\refl{P}$,
and that proves the theorem: i.e. 
\[
|P|^2 - |\refl{P}|^2 = (1-|z|^2)\sum_{j=0}^{n-1} |A_j|^2 +
 (1-|w|^2)\sum_{k=0}^{m-1}|B_k|^2
\]
for some polynomials $A_j$ of degree $(n-1,m)$ and $B_k$ of degree
$(n,m-1)$.
\end{proof}

\end{section}

\begin{section}{Equivalences}\label{sec:equiv}
For those interested in referring back to the original paper of
Geronimo and Woerdeman \cite{Ger-Woe} we have included this section
devoted to demonstrating equivalences between their theorem statements
and ours.  As usual, we let $\rho$ be a probability measure on the
2-torus.

Let $\Lambda_{+} = \{0,\dots,n\}\times \{0,\dots, m\}$ and define the
  doubly Toeplitz matrix corresponding to $\rho$ by 
\[
c^{\rho}_{u-v} = \ip{z^{u_1}w^{u_2}}{z^{v_1}w^{v_2}}_{\rho}
\]

where $u=(u_1,u_2)$ and $v=(v_1,v_2)$ are elements of $\Lambda_{+}$.

First, let us explicitly mention that our nondegeneracy condition is
the same as the one in \cite{Ger-Woe}.

\begin{propose} The following are equivalent.

\begin{enumerate}
\item $\rho$ is nondegenerate (at the $(n,m)$ level).

\item The doubly Toeplitz matrix corresponding to $\rho$ is positive definite:
\[
(c^{\rho}_{u-v})_{u,v \in \Lambda_{+}} > 0.
\]

\end{enumerate}
\end{propose}

\begin{proof} Both conditions are just another way of saying that for
  any complex numbers $a_{jk}$, $j=0,\dots,n$, $k=0,\dots,m$
\[
\ip{\sum_{j=0}^{n} \sum_{k=0}^{m} a_{jk} z^j w^k}{\sum_{j=0}^{n}
  \sum_{k=0}^{m} a_{jk} z^j w^k}_{\rho}
>0
\]
as long as some $a_{jk}$ is nonzero.
\end{proof}

Next, our ``automatic orthogonality condition'' is equivalent to a
couple of different conditions on Toeplitz matrices that are used in
\cite{Ger-Woe}.  The automatic orthogonality condition is also used in
\cite{Ger-Woe2}. 

\begin{propose} 

The following are equivalent.
\begin{enumerate}
\item Automatic orthogonality:
\[
\gkboxurperplt_{\rho} = \gkboxuperplt_{\rho}
\]

\item A low rank condition on the doubly Toeplitz matrix:

\[
\text{rank}(c^{\rho}_{u-v})_{\begin{matrix} u\in \{1,\dots,n\} \times
    \{0,\dots, m\} \\ v \in \{0,\dots, n\}\times \{1,\dots, m\} \end{matrix}} = nm
\]

\item The vanishing of certain blocks of the inverse of the doubly
  Toeplitz matrix:
\[
\left[(c^{\rho}_{u-v})_{u,v \in \Lambda_{+}\setminus \{(0,0)\}}
    \right]^{-1}_{ \begin{matrix} \{1,\dots,n\}\times \{0 \} \\ \{0\}
      \times \{1,\dots, m\} \end{matrix}} =0
\]

\end{enumerate}

\end{propose}

\begin{proof}
The equivalence of the second and third statements follows from
Theorem 2.4.1 in \cite{Ger-Woe}.  We shall prove (1) is equivalent to
(3).  Let $b_{u,v}$ be the $(u,v)$ entry of the inverse of the matrix
\[
(c^{\rho}_{u-v})_{u,v \in \Lambda_{+}\setminus \{(0,0)\}}.
\]
\emph{Warning:} the rows and columns of $(b_{u,v})$ and
$(c^{\rho}_{u-v})$ are each indexed by the set $\Lambda_{+}\setminus
\{(0,0)\}$.  For each $u \in \Lambda_{+}\setminus\{(0,0)\}$ define
\[
\Phi_u \kveczw  := \sum_{v \in \Lambda_{+}\setminus\{(0,0)\}} b_{u,v}
z^{v_1} w^{v_2}
\]
where we write $v=(v_1,v_2)$.  The set $\{\Phi_u\}_{u\in
\Lambda_{+}\setminus\{(0,0)\}}$ is a dual basis for the monomials in
$\gkboxll$, because for $t=(t_1,t_2) \in
\Lambda_{+}\setminus\{(0,0)\}$
\[
\ip{\Phi_u}{z^{t_1} w^{t_2}}_{\rho} = \sum_{v \in
  \Lambda_{+}\setminus\{(0,0)\}} b_{u,v} c^{\rho}_{v-t} =
  \delta_{u,t}.
\]
Therefore, the set $\{\Phi_{(j,0)}\}_{j=1,\dots,n}$ forms a basis for
the subspace $\gkboxllperpup$.  The condition $b_{(j,0),(0,k)} = 0$ for
$j=1,\dots, n$ and $k=1,\dots, m$ is equivalent to saying each
$\Phi_{(j,0)} \in \gkboxlperpup$ and since these polynomials are
linearly independent, this is equivalent to saying $\gkboxllperpup =
\gkboxlperpup$.  By reflection, this is equivalent to $\gkboxrperpdn =
\gkboxurperpdn$.  It is not difficult to see that in general we have
\[
\gkboxrperpdn \oplus_\rho \gkboxurperplt = \gkboxuperplt \oplus_\rho
\gkboxurperpdn
\]
and so $\gkboxrperpdn = \gkboxurperpdn$ if and only if $\gkboxuperplt
= \gkboxurperplt$ by taking orthogonal complements.  Hence, conditions
(1) and (3) are equivalent.  

We conclude by remarking that the above proof shows that the
conditions obtained by switching the roles of $z$ and $w$ are
equivalent to the conditions (1)-(3).
\end{proof}

The following is probably too easy to be a proposition, but we state
it for emphasis.

\begin{propose} Let $\rho$ and $\sigma$ be two probability measures on
  $\mathbb{T}^2$.  The following are equivalent.

\begin{enumerate}
\item The Toeplitz matrices agree: 
\[
c^{\rho}_{u-v} = c^{\sigma}_{u-v} \text{ for } u,v \in \Lambda_{+}.
\]

\item The inner products agree on $\gkbox$:
\[
\ip{P}{Q}_\rho = \ip{P}{Q}_\sigma \text{ for } P,Q \in \gkbox
\]

\item The reproducing kernels agree:
\[
K_\rho \gkbox \kmatzw = K_\sigma \gkbox \kmatzw. 
\]
\end{enumerate}
\end{propose}

\begin{proof} The first two equivalences are easy.  It is also clear
  that (2) implies (3).  If the reproducing kernels agree then the
  inner products agree on linear combinations of kernel functions.
  Since the kernel functions span $\gkbox$, it follows that the inner
  products must agree on all of $\gkbox$.
\end{proof}

With all of this out of the way, let us state the important theorem
(namely Theorem \ref{thm:Ger-Woe2}) in this paper in its original
language.  As mentioned earlier, this theorem is a slight weakening of
Theorem 1.1.2 in \cite{Ger-Woe}.

\begin{theorem}[Geronimo-Woerdeman]
Given is a probability measure $\rho$ on $\mathbb{T}^2$ with
\[
(c^{\rho}_{u-v})_{u,v \in \Lambda_{+}} > 0
\]
 
A unit norm polynomial $p$ in $\gkboxllperp$ is stable and the measure
on $\mathbb{T}^2$
\begin{equation*}
d\mu := \frac{1}{(2\pi i)^2|p\kveczw|^2} \gkmeas
\end{equation*}
has Fourier coefficients $c^{\mu}_u = c^{\rho}_u$ for $u \in
\{-n,\dots,n\} \times \{-m,\dots, m\}$ if and
only if
\[
\text{rank}(c^{\rho}_{u-v})_{\begin{matrix} u\in \{1,\dots,n\} \times
    \{0,\dots, m\} \\ v \in \{0,\dots, n\}\times \{1,\dots, m\} \end{matrix}} = nm.
\]
\end{theorem}

\end{section}

\begin{section}{Conclusions and Commentary}\label{sec:Con}
\begin{subsection}{Three Variables?}
There are several proofs of And\^{o}'s inequality (see \cite{Ando} and
\cite{Ball} for instance).  The work of Geronimo and Woerdeman and
this paper provide two more.  Such a simple looking inequality has
garnered our interest partially because of its connection to bounded
analytic functions and interpolation.  It is also interesting because
the analogous statement for three commuting contractions fails. This
phenomenon is not well understood.  We hope that the notation and
approach of this paper can play some small part in addressing this
problem.  Admittedly, our notation has certain drawbacks for
generalizations to three variables, but we do not think they are
insurmountable.
\end{subsection}

\begin{subsection}{Formulas for reproducing kernels}
It would be interesting to have explicit formulas for the reproducing
kernels involved in the two variable Christoffel-Darboux formula in
the case of Bernstein-Szeg\H{o} measures.  This might allow for a
slicker proof of the Cole-Wermer theorem (i.e. passing to subsequences
might be avoided).  It would also be interesting to study the
probability measures $\mu_r$ in our proof of the Cole-Wermer theorem.
\end{subsection}

\begin{subsection}{Proof of the stability result}
Let us conclude the paper with a proof of the stability result from
the introduction.

\begin{proof}[Proof of Theorem \ref{stablethm}]
The reverse implication is given by Lemma \ref{lemma:zeros}.  To prove
the forward implication, define a Bernstein-Szeg\H{o} measure using
$q$.  The constant $c$ is defined by
\[
\frac{1}{c^2} = \frac{1}{(2 \pi i)^2} \int_{\mathbb{T}^2}
\frac{1}{|q\kveczw|^2} \gkmeas.
\]
By Theorem \ref{bs:thm}, Corollary \ref{cor:ineq}, and Proposition
\ref{bsqprop}, it follows that
\begin{equation*}
|q\kveczw|^2 - |\refl{q}\kveczw|^2 \geq c^2(1-|z|^2)(1-|w|^2).
\end{equation*}
\end{proof}
\end{subsection}
\end{section}

\section*{Symbols and acknowledgments}
The special symbols such as $\gkboxllperp$ used in this article will
be made available on the author's website.  Thanks to Dror Bar-Natan
for making his method of producing new \LaTeX\, symbols available on
his website.  The author would also like to thank John
E. M$^\text{c}$Carthy for his advice at all stages of this research.

\bibliographystyle{plain}
\bibliography{Knese-BernSzeg}











\end{document}